\documentclass{arxiv}

\usepackage{blkarray}
\usepackage{text/kbordermatrix}

\mathtoolsset{showonlyrefs}

\usetikzlibrary{external}
\newcommand{\lagrange}[1]{\mathcal{L}\mleft(#1\mright)}

\newcommand{\pareto}[1]{#1^\star}
\newcommand{\obj}[2][]{f_{#1}\mleft(#2\mright)}
\newcommand{\eq}[2][]{g_{#1}\mleft(#2\mright)}

\newcommand{\lagrangeeq}{\lambda}

\newcommand{\var}{\vc{x}}
\newcommand{\weight}{\vc{w}}

\title{Computing the Pareto Front by Polynomial Elimination, With an Application From System Identification\thanks{This work was supported in part by the KU Leuven Research Fund (grants IBOF/23/064, C3/20/117, and C3I/21/00316); in part by the Fonds Wetenschappelijk Onderzoek (grants S005319 and T001919N); in part by the Departement Economie, Wetenschap \& Innovatie via the Flanders AI Research Program; in part by the Vlaams Agentschap Innoveren \& Ondernemen (grant HBC/2021/0076); and in part by the European Research Council (grant 885682). The work of Christof Vermeersch was supported in part by a FWO PHC Tournesol mobility grant (grant VS00326N).}}

\author[$\ddagger$]{Hans van Rooij}
\author[$\dagger$, $\ddagger$]{Christof Vermeersch}
\author[$\ddagger$]{Marie Deferme}
\author[$\ddagger$]{Bart De Moor}

\affil[$\dagger$]{Corresponding author (\url{christof.vermeersch@esat.kuleuven.be})}
\affil[$\ddagger$]{Center for Dynamical Systems, Signal Processing, and Data Analytics (STADIUS), Dept. of Electrical Engineering (ESAT), KU Leuven, Kasteelpark Arenberg 10, 3001 Leuven, Belgium}

\begin{document}
    \maketitle

    \begin{abstract}
    We propose a novel numerical approach to compute the Pareto front in multivariate polynomial multi-objective optimization problems.
    When the objective functions and (equality) constraints are multivariate polynomials, the Pareto front, which describes the efficient points of the multiple (often conflicting) objective functions, can be interpreted as a subset of a positive-dimensional algebraic variety. 
    By combining the objective functions with weights and considering the weights as additional decision variables, we can eliminate all variables except the objective values and obtain one (or multiple) polynomial equation(s) that describes the Pareto front.
    Unlike sampling-based methods that approximate the Pareto front point-wise, our elimination-based approach yields an explicit algebraic relation between the objective values, representing the Pareto front as a geometric object in the objective space without requiring a predetermined number of sample points.
    Besides numerical examples illustrating the elimination-based approach, we use elimination on a challenging application that originates from system identification, in which we analyze the trade-off between misfit and latency terms when determining the optimal model parameters from measured data.
\end{abstract} 

\section{Introduction}
    \label{sec:introduction}

    Due to possibly conflicting and incommensurable objective functions in multi-objective (or vector) optimization problems, it is generally not possible to find a single solution that is optimal for all objectives simultaneously~\cite{miettinen1999nonlinear}.
    In that sense, such problems do not admit a unique optimum but rather a set of efficient solutions called \emph{Pareto-efficient points}, for which no objective can be improved without worsening another.\footnote{The concept is attributed to both Francis Ysidro Edgeworth and Vilfredo Pareto, two economists from the late 19th, early 20th century. Although Edgeworth presented the idea first, it is named after Pareto, who introduced the notion of economic efficiency in the context of income distribution and welfare economics~\cite{pareto1896cours}. His key insight was that a society's resources are optimally allocated when no individual can be made better off without making another worse off---a condition now known as Pareto efficiency.}
    The corresponding set of objective values, the \emph{Pareto front}, represents the trade-offs inherent to the optimization problem and is, therefore, central in the decision making process.

    When the objective functions and (equality) constraints are multivariate polynomials, the Pareto front can be interpreted as a subset of a positive-dimensional algebraic variety.
    We present in this paper a novel approach to compute in that case the Pareto front by (numerical) polynomial elimination: eliminating all variables except the objective values yields polynomial equations(s) of which a part of the solution set corresponds to the Pareto front.
    
    Existing methods for computing (parts of) the Pareto front typically encode the relative importance of the objective functions before (\textit{a priori}) or after (\textit{a posteriori}) optimization, for instance through weights, constraints, or ordering on the objective functions.
    A survey of such methods is given in~\cite{marler2004survey}.
    A common approach is scalarization, in which the multi-objective optimization problem is reformulated as a family of scalar problems that sample the Pareto front.
    The (convex) weighted sum method is widely used in scalarization because of its simplicity: a weighted combination of the objective functions is minimized for varying weights and each minimizer corresponds to a distinct point on the Pareto front.
    However, it has two well-known limitations: (i) a uniform spread of weights does not generally produce a uniform spread on the Pareto front and (ii) non-convex regions of the Pareto front cannot be detected.
    Moreover, obtaining a finer approximation of the Pareto front increases the computational cost proportionally.

    In contrast, our elimination-based approach treats the weights as additional variables rather than fixed parameters.
    The first-order necessary optimality conditions for the corresponding scalarized, but parametrized, multivariate polynomial multi-objective optimization problem results in a system of multivariate polynomial equations.
    Unlike sampling-based methods that approximate the Pareto front point-wise, eliminating all variables except the objective values yields an explicit algebraic relation between the objective values; the Pareto front is a subset of this implicitly represented geometric object in the objective space.  
    It is no longer necessary to decide in advance on the number of sampling points.  
    Although our approach is still computationally intensive, it does not suffer from the above-mentioned two limitations.
    
    Multi-objective trade-offs arise naturally in many engineering applications (see~\cite{miettinen1999nonlinear, stadler1988multicriteria} and references therein).
    For example, in system identification, determining model parameters that best describe measured data leads inherently to competing objectives.
    Here, we apply our elimination-based approach to the misfit-versus-latency model class, where the mismatch between measured and model-compliant data, which result from external influences such as measured noise, observed disturbances, and unmodeled effects, are represented by a misfit term and an unknown latent input.
    The optimal model parameters are obtained by distributing this mismatch between both terms, resulting in a multi-objective optimization problem that minimizes them simultaneously.
    Although only a limited number of data points is considered, the accompanying numerical example provides both a challenging elimination benchmark (with ten decision variables) and an illustration of the application potential of Pareto fronts in system identification.
    
    \paragraph*{Outline and contributions}

        The remainder of this paper is organized as follows. 
        Firstly, in \cref{sec:paretofronts}, we explain how Pareto fronts appear in multi-objective minimization. 
        We state the problem mathematically and give a motivational example.
        In \cref{sec:methodology}, we present our elimination-based approach for computing an algebraic description of the Pareto front; \emph{this constitutes the main contribution of the paper}. 
        We consider in \cref{sec:numericalexamples} three numerical examples.
        Afterwards, in \cref{sec:application}, we \emph{apply elimination to parameter estimation in the misfit-versus-latency model class, which forms a secondary contribution}. 
        Finally, \cref{sec:conclusion} concludes the paper and outlines directions for future work.

\section{Pareto Fronts in Minimization}
    \label{sec:paretofronts}

    Given $m$ objective functions $\obj[1]{\var}, \obj[2]{\var}, \ldots, \obj[m]{\var}$ in $n$ decision variables $\var = (x_1, x_2, \allowbreak \ldots, x_n) \in \Rset^n$ with $m \geq 2$ and $\obj[i]{\var}: \Rset^n \to \Rset$, the multi-objective minimization problem denotes the simultaneous minimization of each of these objective functions on a given real-valued\footnote{In engineering applications, the decision variables are typically real-valued. For an introductory text that discusses multivariate polynomial minimization with complex-valued decision variables, see~\cite{vermeersch2023multivariate}.} feasible decision space $\mt{X} \subseteq \Rset^n$,
    \begin{equation}
        \label{eq:multiobjectiveminimization}
        \begin{aligned}
            \min_{\var} & \begin{bmatrix} \obj[1]{\var} & \obj[2]{\var} & \ldots & \obj[m]{\var}\end{bmatrix}^\tr, \\
            \sub & \var \in \mt{X}
        \end{aligned}
    \end{equation}
    which is defined as finding all Pareto-efficient points $\pareto{\var}$ in $\mt{X}$.
    The concept of Pareto efficiency (other terms are \emph{admissibility}, \emph{non-inferiority}, and \emph{Pareto optimality}) is used in multi-objective minimization because typically no single global minimum exists, and it is defined as follows.

    \begin{definition}
        A point $\pareto{\var} \in \mt{X}$ is \emph{Pareto-efficient} if there exists no feasible $\var \in \mt{X}$ such that $\obj[i]{\var} \leq \obj[i]{\pareto{\var}}$ for all $i = 1, 2, \ldots, m,$ and $\obj[j]{\var} < \obj[j]{\pareto{\var}}$ for at least one $j$.
    \end{definition}

    The set of objective values of all (possibly infinitely many) Pareto-efficient points is the so-called Pareto front.
    A stricter definition exists that excludes certain degenerate Pareto fronts.
    
    \begin{definition}
        A point $\pareto{\var} \in \mt{X}$ is \emph{properly Pareto-efficient} if it is Pareto-efficient and if there exists a scalar $R > 0$ such that, for each $i$, we have that 
        \begin{equation}
            \label{eq:ratio}
            \frac{\obj[i]{\var} - \obj[i]{\pareto{\var}}}{\obj[j]{\pareto{\var}} - \obj[j]{\var}} \leq R,
        \end{equation}
        for some $j$ such that $\obj[j]{\var} > \obj[j]{\pareto{\var}}$ whenever $\var \in \mt{X}$ and $\obj[i]{\var} < \obj[i]{\pareto{\var}}$.
    \end{definition}
    
    The ratio in~\eqref{eq:ratio} measures the gain in objective function $\obj[i]{\var}$ relative to the loss in objective function $\obj[j]{\var}$.
    Proper Pareto efficiency requires this ratio to be bounded by the same $R$.
    If no finite $R$ exists, which means that there are sequences of feasible points where this ratio grows unbounded, then $\pareto{\var}$ is Pareto-efficient but improper.

    The (convex) weighted sum method rephrases the multi-objective minimization problem~\eqref{eq:multiobjectiveminimization} into a scalar minimization problem given by
    \begin{equation}
        \label{eq:weightedminimization}
        \begin{aligned}
            \min_{\var} & \sum_{i = 1}^m w_i \obj[i]{\var}, \\
            \sub & \var \in \mt{X},
        \end{aligned}
    \end{equation}
    where the weights $w_i$, for $i = 1, 2, \ldots, m$, satisfy the following constraints:
    \begin{equation}
        w_i \ge 0 \eqand \sum_{i=1}^m w_i = 1.
    \end{equation}
    For a specific choice of weights $\weight$, the solution of the weighted minimization problem~\eqref{eq:weightedminimization} gives one point on the Pareto front.
    The proof of the following important theorem is given, for maximization, in~\cite{geoffrion1968proper}.

    \begin{theorem}
        Let $w_i > 0$ be fixed for $i = 1, 2, \ldots m$ with $\sum_{i = 1}^m w_i = 1$.
        If $\pareto{\var}$ is a minimum of the (convex) weighted minimization problem~\eqref{eq:weightedminimization}, then $\pareto{\var}$ is properly Pareto-efficient in~\eqref{eq:multiobjectiveminimization}.
    \end{theorem}
    
    The multi-objective minimization problem potentially has $n_g$ equality constraints $\eq[k]{\var}: \Rset^n \to \Rset$:
    \begin{equation}
        \label{eq:eq}
        \eq[k]{\var} = 0, \quad k = 1, 2, \ldots, n_g, 
    \end{equation}
    which determine the feasible decision space $\mt{X} \subseteq \Rset^n$.

    \begin{definition}
        The \emph{feasible decision space} $\mt{X} \subseteq \Rset^n$, also called \emph{constraint set}, is defined as the set $\left\lbrace \var : \eq[k]{\var} = 0, k = 1, 2, \ldots, n_g \right\rbrace$.
    \end{definition}

    The critical points of the weighted minimization problem~\eqref{eq:weightedminimization} can in that case be found by solving the Karush--Kuhn--Tucker (KKT) conditions~\cite{nocedal2006numerical} for a specific choice of weights $\weight$:
    \begin{equation}
        \label{eq:kktsystem}
        \begin{aligned}
            \nabla_{\var} \lagrange{\var, \vc{\lagrangeeq}} &= \vc{0},\\
            \eq[k]{\var} &= 0, \quad k = 1, 2, \ldots, n_g,
        \end{aligned}
    \end{equation}
    where the Lagrangian $\lagrange{\var,  \vc{\lagrangeeq}}$ is given by 
    \begin{equation}
        \lagrange{\var, \vc{\lagrangeeq}} = \sum_{i = 1}^m w_i \obj[i]{\var} - \sum_{k = 1}^{n_g} \lagrangeeq_k \eq[k]{\var}.
    \end{equation}
    One of the critical points obtained by solving the KKT system~\eqref{eq:kktsystem}, specifically the global minimum, is a Pareto-efficient point of the multi-objective minimization problem~\eqref{eq:multiobjectiveminimization}.

    Since the focus in the paper is on the subclass of multivariate polynomial multi-objective minimization with equality constraints, we restrict ourselves to objective functions and constraints that are multivariate polynomials:
    \begin{equation}
        \obj[i]{\var}, \eq[k]{\var} \in \mathbb{R}[\var].
    \end{equation}
    The Pareto front can then be described by multivariate polynomial equations and is, therefore, contained in a (positive-dimensional) algebraic variety.
    Note that when polynomial inequalities constraints are present, it is a subset of a semi-algebraic set, which itself is contained in a variety.

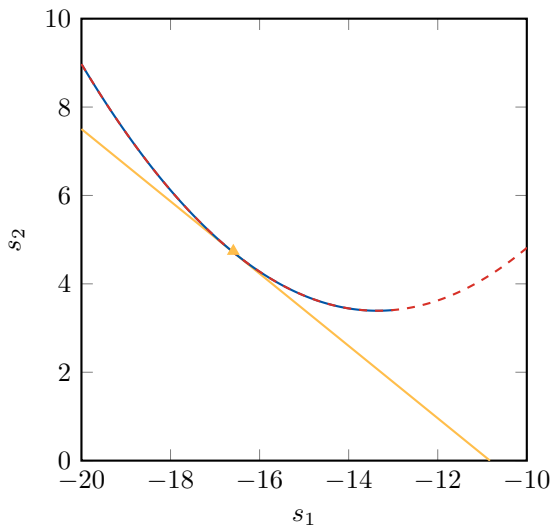
\begin{figure}
        \centering
        \begin{tikzpicture}[trim axis left, trim axis right]
    \begin{axis}[
        figurestyle,
        xmin = -20,
        xmax = -10,
        ymin = 0,
        ymax = 10,
        xlabel = {$s_1$},
        ylabel = {$s_2$},
        ]

        \addplot[yellowline, raw gnuplot] gnuplot {
            set contour base;
            set cntrparam levels discrete 0.001;
            unset surface;
            set view map;
            set isosamples 500;
            set samples 500;
            splot [-20:-10][0:25] -0.031221930729661912017007703579186*x - 0.038170859319640937946899300198466*y - 0.33704195762999307449385511943796;
        };
        \label{plot:toyproblem:tangent}
        

        \addplot[blueline, raw gnuplot] gnuplot {
            set contour base;
            set cntrparam levels discrete 0.001;
            unset surface;
            set view map;
            set isosamples 500;
            set samples 500;
            splot [-20:-1055/79][0:25] 0.0000000000000028033131371785202645696699619293*x^4 + 0.000000000000023290397388464612049574498087168*x^3*y + 0.0000000000020632783676033383812864485662431*x^3 + 0.00000000000000026801477703841669608664233237505*x^2*y^2 - 0.0000000000015137986350555188153066410450265*x^2*y + 0.0048247967210469566365360627457903*x^2 - 0.0000000000000000056378512969246230568387545645237*x*y^3 - 0.00000000000012109627536838196704138681525365*x*y^2 - 0.000000000064048654536486351873669775613962*x*y + 0.12886482279953881402434490155429*x + 0.00000000000000000021889298958106412287115207040116*y^4 - 0.00000000000000054513685232571162941894726827741*y^3 + 0.000000000020072134746814502756495635082956*y^2 - 0.038170860069221428101382542763531*y + 0.99091551820678436257594512426294;
        };
        \label{plot:toyproblem:pareto}




        \addplot[redline, dashed, raw gnuplot] gnuplot {
            set contour base;
            set cntrparam levels discrete 0.001;
            unset surface;
            set view map;
            set isosamples 500;
            set samples 500;
            splot [-20:0][0:10] 0.0000000000000028033131371785202645696699619293*x^4 + 0.000000000000023290397388464612049574498087168*x^3*y + 0.0000000000020632783676033383812864485662431*x^3 + 0.00000000000000026801477703841669608664233237505*x^2*y^2 - 0.0000000000015137986350555188153066410450265*x^2*y + 0.0048247967210469566365360627457903*x^2 - 0.0000000000000000056378512969246230568387545645237*x*y^3 - 0.00000000000012109627536838196704138681525365*x*y^2 - 0.000000000064048654536486351873669775613962*x*y + 0.12886482279953881402434490155429*x + 0.00000000000000000021889298958106412287115207040116*y^4 - 0.00000000000000054513685232571162941894726827741*y^3 + 0.000000000020072134746814502756495635082956*y^2 - 0.038170860069221428101382542763531*y + 0.99091551820678436257594512426294;
        };
        \label{plot:toyproblem:eliminant}

        \addplot [yellowmark, only marks] coordinates{
            (-16.59, 4.74)
        };
        \label{plot:toyproblem:point}

    \end{axis}
\end{tikzpicture}
        \caption{Visualization of the Pareto front for the portfolio optimization problem in \cref{ex:motivationalexample}. The Pareto front (\ref{plot:toyproblem:pareto}) is obtained by considering the zero set of the eliminant polynomial (\ref{plot:toyproblem:eliminant}), which is computed via our novel elimination-based approach in \cref{sec:motivationalexample}. Only the points that correspond to nonnegative weights are contained in the Pareto front. A linear level curve (\ref{plot:toyproblem:tangent}) allows to recover the weight vector and decision variables that correspond to a specific point on the Pareto front (~\ref{plot:toyproblem:point}~).}
        \label{fig:toyproblem}
    \end{figure}

    \begin{example}[Motivational example]      
        \label{ex:motivationalexample}
        Consider the problem of creating an investment portfolio\footnote{As Pareto fronts have their roots in economics, we introduce a small motivational example from that setting, which we use throughout the paper to illustrate the key concepts and our novel approach. The example is based on the work of Nobel prize winner Harry Markowitz~\cite{markowitz1952portfolio}.}, in which one has to decide the amount of money to invest in three different stock options $x_i$, $i = 1, 2, 3$.
        The two conflicting objectives are to maximize revenue, $\obj[1]{\var} = \vc{a}^\tr \var$, and minimize risk, $\obj[2]{\var} = \var^\tr \mt{B} \var$.
        Under the constraint that one has a limited amount of money to invest, $\eq{\var} = \sum_{i = 1}^3 x_i - 100 = 0$, the multi-objective minimization problem corresponds to 
        \begin{align}
            \min_{\var} & \begin{bmatrix} -\vc{a}^\tr \var \\ \var^\tr \mt{B} \var \end{bmatrix} \\
            \sub & \sum_{i = 1}^3 x_i - 100 = 0,
        \end{align}
        where the vector $\vc{a} \in \Rset^3$ contains the expected revenue per stock and the covariance matrix $\mt{B}$ describes the risk associated with every stock.
        Depending on the values of $\vc{a}$ and $\mt{B}$, the objective functions can be conflicting, which is the case for 
        \begin{equation}
            \vc{a} = 
            \begin{bmatrix}
                10 \\ 20 \\ 15
            \end{bmatrix} \cdot 10^{-2} 
            \eqand
            \mt{B} = 
            \begin{bmatrix}
                5 & 1 & 2\\ 1 & 10 & 3\\ 2 & 3 & 7
            \end{bmatrix} \cdot 10^{-4}.
        \end{equation}
        \Cref{fig:toyproblem} displays the Pareto front for these values of $\vc{a}$ and $\mt{B}$ as a solid blue line (\ref{plot:toyproblem:pareto}). 
        The objective variables $s_1$ and $s_2$ denote the negated expected revenue $-\vc{a}^\tr \var$ and expected risk $\var^\tr \mt{B} \var$, respectively.
    \end{example}

\section{Computational Elimination Approach}
    \label{sec:methodology}
    
    We now present our elimination-based approach for computing the Pareto front of multivariate polynomial multi-objective minimization problems. 
    We first describe how the Pareto front can be characterized through a system of multivariate polynomial equations derived from the KKT conditions, and how its computation can be formulated as an elimination problem (\cref{sec:eliminationapproach}). 
    Next (in \cref{sec:numericalalgorithm}), we introduce a numerical elimination algorithm that uses the Macaulay matrix to compute polynomial relations between the objective variables. 
    The interpretation of the resulting eliminant system and corresponding Pareto front are discussed afterwards (\cref{sec:interpretation}). 
    Finally, the complete procedure is illustrated on the motivational example (\cref{sec:motivationalexample}).

    \subsection{Reformulation of the Pareto front as an eliminant}
        \label{sec:eliminationapproach}
        
        If each $\obj[i]{\var},\eq[k]{\var}$ is a multivariate polynomial, then the KKT conditions~\eqref{eq:kktsystem} yield a multivariate polynomial system in the variables $\left(\var, \weight, \vc{\lagrangeeq}\right)$. 
        To obtain a description of the Pareto front in the objective space, we introduce $m$ additional variables
        \begin{equation}
            \vc{s} = \left( s_1, s_2, \ldots, s_m \right),
        \end{equation}
        representing the objective values, and impose the relations
        \begin{equation}
            \label{eq:cost_value_general}
            s_i - \obj[i]{\var} = 0, \quad i = 1, 2, \ldots, m.
        \end{equation}
        Combining the KKT conditions~\eqref{eq:kktsystem} and objective relations~\eqref{eq:cost_value_general} yields the so-called Pareto front (PF) system:
        \begin{equation} 
            \label{eq:pfsystem}
            \begin{aligned}
                \nabla_{\var} \lagrange{\var, \vc{\lagrangeeq}} &= \vc{0}, \\
                \eq[k]{\var} &= 0, \quad k = 1, 2, \ldots, n_g, \\
                s_i - \obj[i]{\var} &= 0, \quad i = 1, 2, \ldots, m.
            \end{aligned}
        \end{equation}
        Eliminating $n' = n + m + n_g$ variables $(\var, \weight, \vc{\lagrangeeq})$ from~\eqref{eq:pfsystem} yields algebraic relations in the objective values $\vc{s}$.
        Hence, computing the Pareto front can be formulated as an elimination problem. 
        The \emph{eliminant system}, of which the solutions contain the Pareto front, may consist of a single polynomial equation or a multivariate polynomial system.

        The elimination of $(\var, \weight, \vc{\lagrangeeq})$ can be done symbolically or numerically.
        Numerical elimination is a relatively recent development compared to the long-standing symbolic tradition. 
        In contrast to symbolic approaches, numerical formulations rely on techniques from numerical linear algebra at every stage of the computation, allowing floating-point arithmetic~\cite{batselier2013geometrical, zeng2008numerical, vanrooij2026numerical}.

    \subsection{Numerical elimination algorithm}
        \label{sec:numericalalgorithm}
    
        Our numerical elimination algorithm uses Macaulay matrices $\mt{M}_{d} \in \mathbb{R}^{p_d \times q_d}$ of increasing degree $d$, which are built from the coefficients of the PF system~\eqref{eq:pfsystem} and monomial multiples up to total degree $d$~\cite{vermeersch2023block}. 
        To construct a Macaulay matrix, we take each polynomial in the PF system~\eqref{eq:pfsystem} and multiply it by all monomials of increasing total degree so that the highest total degree does not exceed $d$, and we collect the coefficients of these polynomials in a structured matrix with columns labeled by the monomials~\cite{vermeersch2023block}.
        Its dimensions $p_d \times q_d$ depend on the degree $d$.
        
        Each row of the Macaulay matrix $\mt{M}_d$ represents one polynomial (either an original equation or a monomial multiple) expressed in a fixed monomial basis~\cite{vermeersch2023block}. 
        Hence, the row space of $\mt{M}_d$, denoted by $\mcl{R}_d$, consists of polynomial combinations of the polynomials, $r_i(\vc{z})$ with total degree $d_i$, in the PF system~\eqref{eq:pfsystem}:
        \begin{equation}
            \mcl{R}_d = \left\{\sum_{i = 1}^{n'} c_i(\vc{z}) \, r_i(\vc{z}) : c_i(\vc{z}) \in \Rset[\vc{z}]_{\le d - d_i} \right\}.
        \end{equation}
        Here, $\Rset[\vc{z}]_{\le d - d_i}$ denotes the set of polynomials in 
        $\vc{z} = \left(\var, \weight, \vc{\lagrangeeq}, \vc{s}\right)$ with real coefficients of degree at most $d - d_i$.
        Thus, the row space $\mcl{R}_d$ can be interpreted as the space of all polynomial combinations of the PF system with degree bounded by $d$.
                    
        The goal of our numerical elimination algorithm is to remove all variables except the objective variables $\vc{s}$. 
        To achieve this goal, we construct polynomial combinations of the original system that no longer contain the elimination variables. 
        Since each vector in the row space $\mcl{R}_d$ represents such a polynomial combination, this amounts to finding vectors whose elements corresponding to monomials involving elimination variables are zero. 
        Equivalently, we seek vectors whose nonzero entries correspond only to monomials in $\vc{s}$. 
        We construct a so-called \emph{zero pattern vector} $\vc{t} \in \Rset^{q_d}$ that ``forces'' coefficients to zero in the row space.
        The idea is illustrated in the following example.
        
        \begin{example}
            \label{ex:zeropattern}
            Let $\var = (x_1, x_2, x_3)$ and consider two objectives ($m = 2$) without constraints.
            The variables are $(x_1, x_2, x_3, w_1, w_2, s_1, s_2)$. 
            We eliminate all decision and weight variables, so that only the objective variables $(s_1, s_2)$ remain.
            The predetermined zero pattern in the row space corresponds to
            \begin{equation}
                \vc{t}^\tr =
                \begin{blockarray}{cccccccc}
                    1 & x_1 & \cdots & w_2 & s_1 & s_2 & x_1^2 & \\
                    \begin{block}{[cccccccc]}
                        \bullet & 0 & 0 & 0 & \bullet & \bullet & 0  & \cdots \\
                    \end{block}
                \end{blockarray},
            \end{equation}
            where $\bullet$ denotes an arbitrary real number. 
            This vector represents a linear combination of the rows of $\mt{M}_d$ in which all coefficients corresponding to monomials containing elimination variables are zero, while coefficients of monomials in only $(s_1, s_2)$ may be nonzero. 
            As a result, the corresponding polynomial(s) depend solely on $(s_1, s_2)$.
        \end{example}
        
        The zero pattern vector $\vc{t}$ can be parametrized as
        \begin{equation}
            \vc{t}^\tr =
            \begin{bmatrix}
                \bullet &\bullet &\bullet & \cdots
            \end{bmatrix}
            \underbrace{\begin{bmatrix}
                1 & 0 & 0 & 0 & 0 & 0 & 0 & 0 & \cdots  \\
                0 & 0 & 0 & 0 & 0 & 0 & 1 & 0 & \cdots  \\
                0 & 0 & 0 & 0 & 0 & 0 & 0 & 1 & \cdots  \\
                \vdots & \vdots & \vdots & \vdots & \vdots & \vdots & \vdots & \vdots & \ddots \\
            \end{bmatrix}}_{\mt{E}_d}.
        \end{equation}
        The rows of $\mt{E}_d$ form an orthonormal basis of the elimination vector space $\mcl{E}_d$.
        The matrix $\mt{E}_d$ acts as a selector that keeps only the entries corresponding to monomials in the decision variables $\vc{s}$ and sets all other entries to zero. 
        Consequently, $\mcl{E}_d$ consists of all coefficient vectors that have this prescribed zero pattern, i.e., vectors whose nonzero entries correspond only to monomials in $\vc{s}$.
        
        Consequently, elimination reduces to finding vectors that lie in both $\mcl{R}_d$ and $\mcl{E}_d$, i.e., to computing their intersection.
        Vectors in $\mcl{R}_d$ correspond to polynomial combinations of the original system, while vectors in $\mcl{E}_d$ encode the requirement that only monomials in the objective variables $\vc{s}$ remain.
        Therefore, vectors in the intersection $\mcl{R}_d \cap \mcl{E}_d$ represent polynomial combinations of the PF system~\eqref{eq:pfsystem} in which all elimination variables vanish, yielding relations purely in $\vc{s}$.
        A nontrivial intersection requires
        \begin{equation}
            \label{eq:dim}
            \dim(\mcl{R}_d \cap \mcl{E}_d) > 0.
        \end{equation}
        Using Grassmann’s dimension theorem~\cite{fearnleysander1979hermann}, the dimension of the intersection is given by
        \begin{equation}
            \dim(\mcl{R}_d \cap \mcl{E}_d) = \rank(\mt{M}_d) - \rank(\mt{N}_d),
        \end{equation}
        where $\mt{N}_d \in \mathbb{R}^{p_d \times q_d'}$ is a submatrix of the Macaulay matrix $\mt{M}_d$ that consists of the $q_d'$ columns associated with at least one variable that need to be eliminated.\footnote{The presented approach uses the left singular vectors of $\mt{N}_d$, which is a submatrix of the Macaulay matrix. Although not explained in this paper, an alternative approach via principal angles and directions between $\mathcal{R}_d$ and $\mathcal{E}_d$ is computationally more efficient~\cite{vanrooij2026numerical}.}
        A positive dimension, i.e.,
        \begin{equation}
            \label{eq:rankcondition}
            \rank (\mt{M}_d) > \rank(\mt{N}_d),
        \end{equation}
        guarantees the existence of an eliminant system, meaning that there exist nontrivial combinations of the original polynomials in which all eliminated variables cancel out. 

        The degree $d$ is increased until the rank condition in~\eqref{eq:rankcondition} is satisfied, since higher degrees include more monomial multiples of the original equations, allowing additional polynomial combinations that can eliminate the elimination variables.
        
        For that degree, a basis $\mt{V} \in \mathbb{R}^{p_d \times l}$ for the left null space of $\mt{N}_d$ is computed, such that
        \begin{equation}
            \mt{V}^\tr \mt{N}_d = \mt{0},
        \end{equation}
        A left multiplication of the Macaulay matrix $\mt{M}_d$ by $\mt{V}^\tr$ thus results in linear combinations of the rows of $\mt{M}_d$, while eliminating all monomials containing one or more of the variables to be removed.
        The result, $\mt{S}^\tr = \mt{V}^\tr \mt{M}_d \in \mathbb{R}^{l \times q_d}$, only involves the objective variables $\vc{s}$.
        
        The eliminant system is then constructed by interpreting the $l'$ nonzero columns of $\mt{S}$ as coefficient vectors of multivariate polynomials
        \begin{equation}
            \label{eq:elsystem}
            t_1(\vc{s}) = \cdots = t_{l'}(\vc{s}) = 0.
        \end{equation}
        These polynomials form the eliminant system, of which the solution set contains the Pareto front.
        
        Computing the dimension of the intersection between $\mcl{R}_d$ and $\mcl{E}_d$ involves two rank computations per degree $d$, which constitute the main computational complexity of the elimination-based approach.
        The first one requires $\mathcal{O}(p_d q_d \min(p_d, q_d))$ operations, while the second one has complexity $\mathcal{O}(p_d q_d' \min(p_d, q_d'))$.
        Once the degree $d$ is found for which~\eqref{eq:rankcondition} is fulfilled, an orthogonal basis matrix $\mt{V}$ for the left null space of $\mt{N}_d$ is computed, which requires again $\mathcal{O}(p_d q_d' \min(p_d, q_d'))$ operations. 
        Finally, a matrix-matrix multiplication, $\mt{S}^\tr = \mt{V}^\tr \mt{M}_d$, yields the vector space that contains the representation of the eliminant system, costing $\mathcal{O}(l p_d q_d)$ operations.

    \subsection{Geometric interpretation of the Pareto front}
        \label{sec:interpretation}
        
        Let $\mathcal{P} \subset \Rset^m$ denote the Pareto front in the objective space, obtained after elimination as (a subset of) the algebraic variety defined by~\eqref{eq:elsystem}.
        The Jacobian matrix of the eliminant system is
        \begin{equation}
            J_t(\pareto{\vc{s}}) =
            \begin{bmatrix}
                \nabla t_1(\pareto{\vc{s}})^\tr \\
                \vdots \\
                \nabla t_{l'}(\pareto{\vc{s}})^\tr
            \end{bmatrix},
        \end{equation}
        and the tangent space of $\mathcal{P}$ at $\pareto{\vc{s}}$ is then given by $\null \left(J_t(\pareto{\vc{s}})\right)$, while the normal space is $\row \left(J_t(\pareto{\vc{s}})\right)$, where $\null$ and $\row$ denote the right null space and row space of $J_t(\pareto{\vc{s}})$, respectively.
        
        At a Pareto-efficient point $\pareto{\vc{s}}$, there exists a vector $\pareto{\weight} \in \mathbb{R}^m$ such that
        \begin{equation}
            (\pareto{\weight})^\tr \vc{s} \ge (\pareto{\weight})^\tr \pareto{\vc{s}},
        \end{equation}
        for all feasible points $\vc{s}$ in the objective space.
        The hyperplane $(\pareto{\weight})^\tr \vc{s} = (\pareto{\weight})^\tr \pareto{\vc{s}}$ passes through $\pareto{\vc{s}}$ and all feasible points lie on one side of it. 
        Geometrically, this hyperplane is tangent to the Pareto front at $\pareto{\vc{s}}$; its normal vector must belong to the normal space of $\mathcal{P}$ at that point. 
        Therefore,
        \begin{equation}
            \pareto{\weight} \in \row \left(J_t(\pareto{\vc{s}})\right).
        \end{equation}
        In the case where the Pareto front is described by a single polynomial equation, $t_1(\vc{s}) = 0$, this reduces to checking whether $\pareto{\weight}$ is parallel to $\nabla t_1(\pareto{\vc{s}})$.
        Alternatively, when the eliminant system~\eqref{eq:elsystem} consist of multiple polynomials, the weight vector lies in the span of their gradients, i.e.,
        \begin{equation}
            \pareto{\weight} = \sum_{i = 1}^{l'} \alpha_i \nabla t_i(\pareto{\vc{s}}), \quad \alpha_i \in \Rset.
        \end{equation}
        
        Since we consider a multi-objective minimization problem, only weight vectors with strictly positive components are of interest.
        Not every point of the algebraic variety thus corresponds to a Pareto-efficient point.
        The Pareto front $\mathcal{P}$ consists of those points $\pareto{\vc{s}}$ for which the normal space $\row \left(J_t(\pareto{\vc{s}})\right)$ contains a vector with only 
        nonnegative components.

        For fixed weights $\pareto{\weight}$, the corresponding decision variables $\pareto{\var}$ can be obtained from the first-order optimality conditions given by the KKT system~\eqref{eq:kktsystem}. 
        Substituting $\pareto{\weight}$ results in a 
        multivariate polynomial system that is generically zero-dimensional. 
        Solving the substituted system yields all critical points $\pareto{\var}$ and the corresponding Lagrange multipliers $\pareto{\vc{\lagrangeeq}}$.

    \subsection{Step-by-step solution of the motivational example}
        \label{sec:motivationalexample}

        \begin{example}
            The PF system associated with \cref{ex:motivationalexample} takes the form
            \begin{equation}
                \label{eq:toyproblemsystem}
                \begin{aligned}
                    -w_1 \vc{a} + 2w_2 \mt{B}\vc{x} + \lambda \vc{1} &= \vc{0},\\
                    x_1 + x_2 + x_3 - 100 &= 0,\\
                    s_1 + \vc{x}^\tr \vc{a} &= 0,\\
                    s_2 - \vc{x}^\tr B \vc{x} &=0,
                \end{aligned}
            \end{equation}
            with $\vc{1} = \begin{bmatrix} 1 & 1 & 1 \end{bmatrix}^\tr$.
            Since we are using convex weights, $w_2$ is replaced by $1 - w_1$. 
            We have a system of six equations in seven variables and eliminate the variables $x_1$, $x_2$, $x_3$, $w_1$, and $\lambda$.
            The elimination-based approach requires a $\num{384} \times \num{330}$ Macaulay matrix of degree four and yields a bivariate polynomial relation between $s_1$ and $s_2$. 
            The zero set of this eliminant polynomial is shown as the red dashed curve (\ref{plot:toyproblem:eliminant}) in \cref{fig:toyproblem}.
            Every point on the true Pareto front is contained in this zero set. 
            The converse does not hold: not all points in the zero set correspond to Pareto-efficient solutions. 
            A discrepancy can arise because the KKT conditions are necessary but not sufficient, and therefore they also capture stationary points that are not minima. 
            However, this is not the case here. 
            The difference arises from omitting the inequality constraints $w_i \geq 0$.

            Suppose that we are interested in the point $\pareto{\vc{s}} = (-16.59, 4.74)$ on the Pareto front, which is indicated by a yellow triangle (~\ref{plot:toyproblem:point}~) on \cref{fig:toyproblem}. 
            The gradient of the eliminant polynomial at this point equals $\begin{bmatrix} -3.12 & -3.82 \end{bmatrix}^\tr \cdot 10^{-2}$, which yields, after normalization, the weights $\pareto{\weight} =  (0.45, 0.55)$. 
            The level curve $(\pareto{\weight})^\tr \vc{s} = (\pareto{\weight})^\tr \pareto{\vc{s}}$ is also shown on \cref{fig:toyproblem} in yellow (\ref{plot:toyproblem:tangent}). 
            As can be seen, this curve is tangent to the Pareto front at $\pareto{\vc{s}}$. 
            Substituting $\pareto{\vc{s}}$ and $\pareto{\weight}$ in the PF system in~\eqref{eq:toyproblemsystem}, and solving the resulting system, gives the corresponding decision variable $\pareto{\var} = (18.18, 50.00, 31.82)$. 
            It can easily be verified that this point indeed corresponds to the objective values in $\pareto{\vc{s}}$.
        \end{example}

\section{Various Numerical Examples}
    \label{sec:numericalexamples}

    To evaluate the practical behavior of the proposed elimination-based approach, we present three numerical examples.\footnote{We use the \macaulaylab toolbox to perform our numerical examples, which is a \matlab implementation of to construct, manipulate, and use Macaulay matrices~\cite{vermeersch2025solving}.}
    Although these examples appear simple, they already illustrate how rapidly the Macaulay matrix (and therefore the computational complexity) grows when the degree $d$ increases.
    For each numerical example, we compute the eliminant system that describes the Pareto front and verify its correctness by substituting points on the sampled Pareto front into the algebraic expression. 
    The characteristics are reported in \cref{table:numericalresults}.

    \begin{example}
        \label{ex:numericalexample1}
        Consider the unconstrained multivariate polynomial tri-objective optimization problem
        \begin{equation}
            \min_{\var} 
            \begin{bmatrix}
                (x_1 - 3)^2 + (x_2-2)^2\\
                x_1 + x_2\\
                x_1 + 2x_2
            \end{bmatrix}.
        \end{equation}
        The first objective function has a unique minimizer $\pareto{\var} = (3, 2)$, while the other two objective functions are linear and favor different directions in $\Rset^2$. 
        As a result, the three objective functions are conflicting, and the Pareto front is a subset of a two-dimensional surface in the $(s_1, s_2, s_3)$-space. 
        The numerical elimination algorithm yields an eliminant polynomial in the objective variables whose zero set describes this surface.
        After rounding the coefficients for readability, the eliminant can be written as $5s_2^2 - 6s_2s_3 + 2s_3^2 - s_1 - 8s_2 + 2s_3 + 13=0$.
    \end{example}

    The next numerical example illustrates how the Macaulay matrix can grow quickly, even when the underlying optimization problem is low-dimensional and structurally simple.

    \begin{example} 
        \label{ex:numericalexample2}
        Consider the constrained multivariate polynomial bi-objective optimization problem
        \begin{align}
            \min_{\var} & 
            \begin{bmatrix}
                -x_1^3 - x_2^3\\
                x_1^2 - x_2^2
            \end{bmatrix}, \\
            \sub & x_1^2 + (x_2 + 1)^2 - 1 = 0.
        \end{align}
        Applying our elimination-based approach to this problem results in a polynomial of degree six: $s_2^6 - 12 s_2^5 - 12 s_1 s_2^4 + 16 s_1^4 + 48 s_1^3 s_2 + 48 s_1^2 s_2^2 + 32 s_1 s_2^3 + 48 s_2^4 - 32 s_1^3 - 48 s_1^2 s_2$. 
       The degree of the Macaulay matrix, however, has to be eight, leading to a $\num{3234} \times \num{3003}$ matrix. 
    \end{example}

    Not every Pareto front is contained in the solution set of a single multivariate polynomial.

    \begin{example}
        \label{ex:numericalexample3}
        Consider the unconstrained multivariate polynomial tri-objective problem
        \begin{align}
            \min_{\var} & 
            \begin{bmatrix}
                x^2 \\
                (x-1)^2 \\
                (x-2)^2
            \end{bmatrix},
        \end{align}
        over a one-dimensional decision space. 
        The Pareto front is a one-dimensional curve embedded in $\Rset^3$. 
        Its algebraic description is given by the two independent polynomial equations
        \begin{align}
            \label{eq:symbolicequations}
            t_1(s_1,s_2,s_3) &= s_1 - 2s_2 + s_3 - 2 = 0, \\
            t_2(s_1,s_2,s_3) &= s_2^2 - 2s_2s_3 + s_3^2 - 2s_2 - 2s_3 + 1 = 0.
        \end{align}
        The numerical elimination algorithm for a Macaulay matrix of degree three does not produce this system exactly, but instead returns five polynomial equations that have the same zero set. 
        The numerical elimination procedure returns a set of linearly independent polynomials, but it does not enforce algebraic independence. 
        As a consequence, the output may contain more equations than strictly necessary, even though they define the same zero set as the symbolically obtained relations in~\eqref{eq:symbolicequations}. 
        The residual in \cref{table:numericalresults} is calculated using the numerically obtained eliminant system.
    \end{example}

    \begin{table}
        \centering
        \caption{Characteristics of the numerical examples.  
        Here, $m$ and $n$ denote the number of objectives and decision variables, respectively; $\#z$ is the total number of variables in the PF system; $d$ is the minimal required degree of the Macaulay matrix; $p_d$ and $q_d$ are the matrix dimensions for that degree; and the final column reports the maximum residual obtained by evaluating the eliminant system at sampled points on the Pareto front.}
        \label{table:numericalresults}
        \begin{tabular}{r|ccccccc}
            \toprule
            problem & $m$ & $n$ & $\# z$ & $d$ & $p_d$ & $q_d$ & $\max \norm{\vc{t}(\pareto{\vc{s}})}_2$ \\
            \midrule 
            \cref{ex:numericalexample1} & $\num{3}$ & $\num{2}$ & $\num{7}$ & $\num{2}$ & $\num{19}$ & $\num{36}$ & $\num{3.16e-15}$ \\
            \cref{ex:numericalexample2} & $\num{2}$ & $\num{2}$ & $\num{6}$ & $\num{8}$ & $\num{3234}$ & $\num{3003}$ & $\num{2.79e-12}$ \\
            \cref{ex:numericalexample3} & $\num{3}$ & $\num{1}$ & $\num{6}$ & $\num{3}$ & $\num{49}$ & $\num{84}$ & $\num{4.73e-15}$ \\
            \bottomrule
        \end{tabular}
    \end{table}
    
\section{System Identification Problem}
    \label{sec:application}

    We apply elimination to the misfit-versus-latency model class, where finding the model parameters that ``best'' describe measured data is inherently a minimization problem with two competing objective functions.
    This example also serves to illustrate the potential of elimination-based techniques for computing Pareto fronts.
    
    We start firstly by giving a short background on the misfit-versus-latency model class (\cref{sec:misfitlatencymodel}).
    Afterwards, we explain how the Pareto front enters the picture (\cref{sec:misfitlatencyparetofront}) and solve the identification problem for some given output data (\cref{sec:misfitlatencyexample}).
    Notice that our approach is at this point not useful for real-life system identification problems with hundreds of data points, due to scale of the multivariate polynomial optimization problem and its corresponding computational complexity.
        
    \subsection{Misfit-versus-latency modeling}
        \label{sec:misfitlatencymodel}

        Misfit-versus-latency modeling~\cite{lemmerling2001misfit} is a system identification framework for (dynamic) linear time-invariant models that combines the so-called ``latency'' approach in which residuals are introduced to explain the measured data (i.e., the prediction error approach~\cite{ljung1999system} containing models such as autoregressive moving-average), and ``misfit'' approach that considers an error term between the measured and model-compliant data (i.e., the behavioral approach~\cite{willems1986time} containing models such as total least squares or errors-in-variables).
        The motivation for this modeling approach lies in the fact that there is no reason to exclude a misfit of the measurements in the presence of a latency term, and vice versa.

        \begin{figure}
            \centering
            \begin{tikzpicture}
    \draw[thick] (0, 0) rectangle (1.5cm, 1cm);
    \node at (0.75cm, 0.5cm) {$\frac{b(q)}{a(q)}$};
    \draw[thick] (0, 1.5cm) rectangle (1.5cm, 2.5cm);
    \node at (0.75cm, 2cm) {$\frac{c(q)}{a(q)}$};

    \draw[thick] (-1cm, -0.5cm) circle (0.25cm);
    \node at (-1cm, -0.5cm) {$+$};
    \draw[thick] (2.5cm, 0.5cm) circle (0.25cm);
    \node at (2.5cm, 0.5cm) {$+$};
    \draw[thick] (3.5cm, -0.5cm) circle (0.25cm);
    \node at (3.5cm, -0.5cm) {$+$};

    \draw[->, thick] (-2cm, 0.5cm) -- (-0, 0.5cm);
    \node[anchor = east] at (-2cm, 0.5cm) {$\vchat{u}\vphantom{\vchat{y}}$};
    \draw[->, thick] (-2cm, 2cm) -- (-0, 2cm);
    \node[anchor = east] at (-2cm, 2cm) {$\vc{e}$};
    \node[font = \itshape] at (-1cm, 1.75cm) {(latency)};
    \draw[->, thick] (1.5cm, 0.5cm) -- (2.25cm, 0.5cm);
    \draw[->, thick] (2.75cm, 0.5cm) -- (4.5cm, 0.5cm);
    \node[anchor = west] at (4.5cm, 0.5cm) {$\vchat{y}$};

    \draw[->, thick] (1.5cm, 2cm) -- (2.5cm, 2cm) -- (2.5cm, 0.75cm);
    
    \draw[->, thick] (-1cm, 0.5cm) -- (-1cm, -0.25cm);
    \draw[->, thick] (3.5cm, 0.5cm) -- (3.5cm, -0.25cm);

    \draw[->, thick] (-1cm, -0.75cm) -- (-1cm, -1.5cm);
    \node[anchor = north] at  (-1cm, -1.5cm) {$\vc{u}$};
    \draw[->, thick] (3.5cm, -0.75cm) -- (3.5cm, -1.5cm);
    \node[anchor = north] at  (3.5cm, -1.5cm) {$\vc{y}$};

    \draw[->, thick] (0, -0.5cm) -- (-0.75cm, -0.5cm);
    \node[anchor = west] at (0, -0.5cm) {$\vctilde{u}\vphantom{\vctilde{y}}$};
    \draw[->, thick] (2.5cm, -0.5cm) -- (3.25cm, -0.5cm);
    \node[anchor = east] at (2.5cm, -0.5cm) {$\vctilde{y}$};
    \node[font = \itshape] at (1.25cm, -0.5cm) {(misfit)};

    \draw[thick, dashed] (-1.5cm, -1.25cm) rectangle (4cm, -2.25cm); 
    \node[anchor = north, font = \itshape] at (1.25cm, -1.5cm) {(measured data)};
    
\end{tikzpicture}
            \caption{Schematic representation of the misfit-latency model class: the measured data ($\vc{u}, \vc{y}$) is described by a latent input term ($\vc{e}$) and misfit term ($\vctilde{u}, \vctilde{y}$), such that the model-compliant data ($\vchat{u}, \vchat{y}$) satisfies a linear time-invariant model equation of the form $a(q) \vchat{y} = b(q) \vchat{u} + c(q) \vc{e}$.}
            \label{fig:diagram}
        \end{figure}
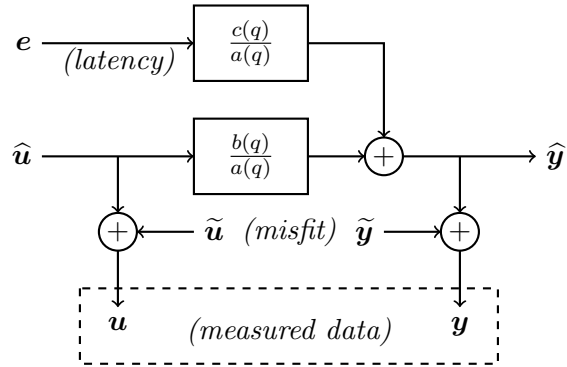

        A schematic representation of the misfit-versus-latency model class is given in \cref{fig:diagram}.
        It considers model equations of the form
        \begin{equation}
            \label{eq:modelequation}
            a(q) \vchat{y} = b(q) \vchat{u} + c(q) \vc{e},
        \end{equation}
        in which $\vchat{u} \in \Rset^N$ and $\vchat{y} \in \Rset^N$ are model-compliant input and output data of length $N$, respectively. 
        The polynomials $a(q)$, $b(q)$, and $c(q)$, respectively of degrees $n_a$, $n_b$ and $n_c$, are polynomials in the delay operator $q$ and constitute the model dynamics.
        The vector $\vc{e} \in \Rset^{N - n_a + n_c}$ represents the unkown, latent input.
        Between the measured input-output data $\vc{u}$, $\vc{y}$ and the model-compliant input-output data $\vchat{u}$, $\vchat{y}$, there is a misfit
        \begin{equation}
            \label{eq:misfitequation}
            \vctilde{y} = \vc{y} - \vchat{y}, \quad \vctilde{u} = \vc{u} - \vchat{u}.
        \end{equation}

        When identifying the model parameters, the goal is to minimize both the misfit and latency introduced in the model, leading to the following objective function
        \begin{equation}
            \label{eq:misfitlatencyobjective}
            \frac{1}{2} \left( \alpha \norm{\vctilde{y}}^2_2 + \beta \norm{\vctilde{u}}^2_2 + \gamma \norm{\vc{e}}^2_2 \right),
        \end{equation}
        with $\norm{\cdot}_2$ denoting the Euclidean norm. 
        The scalars $\alpha$, $\beta$, and $\gamma$ are nonnegative hyperparameters that determine the relative importance of the misfit and latency terms, and they can be infinity as well.
        Several well-known model classes arise for particular choices of these hyperparameters (see~\cite[Table 1]{lemmerling2001misfit}).
        Together with~\eqref{eq:modelequation} and~\eqref{eq:misfitequation}, the objective function~\eqref{eq:misfitlatencyobjective} defines a constrained multivariate polynomial optimization problem.
        Using the KKT conditions, a system of multivariate polynomial equations can be constructed.
        For fixed values of $\alpha, \beta$ and $\gamma$, this multivariate polynomial system can be solved using any polynomial root-finding method.

    \subsection{Relation to the Pareto front}
        \label{sec:misfitlatencyparetofront}

        When considering misfit-versus-latency models, we actually have a trade-off between minimizing the output misfit, input misfit, and latency term. 
        Rather than fixing the hyperparameters $\alpha, \beta$ and $\gamma$ in advance, we can approach misfit-versus-latency modeling as a multi-objective polynomial minimization problem:
        \begin{equation}
            \label{eq:misfitlatencymultiobjectiveproblem}
            \begin{aligned}
                \min_{\var} & \frac{1}{2} \begin{bmatrix} \norm{\vctilde{y}}_2^2 & \norm{\vctilde{u}}_2^2 & \norm{\vc{e}}_2^2 \end{bmatrix}^\tr \\
                \sub & a(q) \vchat{y} = b(q) \vchat{u} + c(q) \vc{e},
            \end{aligned}
        \end{equation}
        where the decision variables $\var$ are $\vc{a}$, $\vc{b}$, $\vc{c}$, $\vc{e}$, $\vchat{y}$, and $\vchat{u}$.
        Via elimination, we aim to characterize the Pareto front that describes the optimal misfit and latency for each possible combination of hyperparameters, which constitute, after normalization, the weights $\weight$.

    \subsection{Applied numerical example}
        \label{sec:misfitlatencyexample}

        As an applied numerical example, we consider a simplified version of the model in~\eqref{eq:modelequation} that is autonomous and has no moving-average term ($n_c = 0$). 
        The minimization problem in~\eqref{eq:misfitlatencymultiobjectiveproblem} reduces to the following weighted minimization problem:
        \begin{equation}
            \label{eq:numericalsystemoptimization}
            \begin{aligned}
                \min_{\vc{a}, \vchat{y}, \vc{e}} &\ \frac{1}{2} \left( \alpha \norm{\vctilde{y}}^2 + (1 - \alpha) \norm{\vc{e}}^2 \right), \\
                \sub &\ \mthat{Y} \vc{a} = \vc{e}
            \end{aligned}
        \end{equation}
        We only retain one hyperparameter and replace the other one ($\gamma$ in this case) by $1 - \alpha$, resulting in a convex combination.
        The matrix $\mthat{Y} \in \Rset^{(N - n_a) \times (n_a + 1)}$ is a Hankel matrix constructed from the model-compliant output data $\vchat{y}$ and $\mthat{Y} \vc{a} = \vc{e}$ is identical to considering $a(q) \vchat{y} = \vc{e}$ for all data points.
        The Lagrangian is equal to 
        \begin{equation}
            \begin{aligned}
                \lagrange{\vc{a}, \vchat{y}, \vc{e}, \vc{\lagrangeeq}}
                &= \frac{1}{2} \left( \alpha \norm{\vctilde{y}}^2_2 + (1 - \alpha) \norm{\vc{e}}^2_2 \right) \\
                &\quad - \vc{\lagrangeeq}^\tr \left( \mthat{Y} \vc{a} - \vc{e} \right),
            \end{aligned}
        \end{equation}
        where $\vc{\lagrangeeq} \in \Rset^{N - n_a}$ are the Lagrange multipliers for the equality constraint. 
        The corresponding KKT system is 
        \begin{equation} 
            \label{eq:misfitlatencykkt}
            \begin{aligned}
                - (\mthat{Y}')^\tr \vc{\lagrangeeq} &= \vc{0}, \\
                \alpha (\vc{y}-\vchat{y}) - \mt{T}_{\vc{a}}^\tr \vc{\lagrangeeq} &= \vc{0}, \\
                (1-\alpha)\vc{e} + \vc{\lagrangeeq} &= \vc{0}, \\
                \mthat{Y}\vc{a} - \vc{e} &= \vc{0},
            \end{aligned}
        \end{equation}
        in which $\mthat{Y}'$ denotes the Hankel matrix $\mthat{Y}$ without its right-most column and $\mt{T}_{\vc{a}}$ is a banded Toeplitz matrix with $\vc{a}$ on its band diagonal.
        The PF system can be constructed from~\eqref{eq:misfitlatencykkt} supplemented with
        \begin{align}
            s_1 - \norm{\vctilde{y}}^2_2 &= 0, \\
            s_2 - \norm{\vc{e}}^2_2 &= 0.
        \end{align}
        Eliminating all decision variables except $s_1$ and $s_2$ from this multivariate polynomial system with $3 N - 2 n_a + 2$ equations yields an algebraic expression that contains the Pareto front of the identification problem.

\begin{figure}
            \centering
            \input{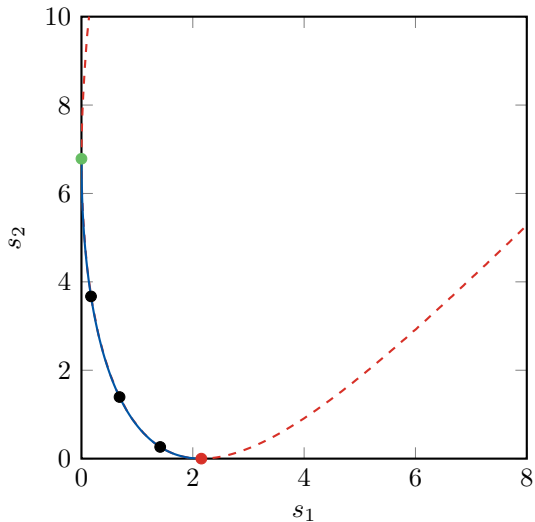}
            \caption{Visualization of the Pareto front for the misfit-versus-latency modeling problem in \cref{ex:mis_vs_lat}. The Pareto front (\ref{plot:misfitlatency:pareto}) is a subset of the eliminant (\ref{plot:misfitlatency:eliminant}) obtained through the elimination approach described in \cref{sec:numericalalgorithm}, and it is validated by sampling Pareto-efficient points (~\ref{plot:misfitlatency:sampled}~). The two intersection points with the axes correspond to a pure autoregressive model without misfit (~\ref{plot:misfitlatency:green}~) and a pure autonomous model with output misfit (~\ref{plot:misfitlatency:red}~), which can also be obtained using dedicated identification algorithms.}
            \label{fig:misfitlatency}
        \end{figure}

        \begin{example}\label{ex:mis_vs_lat}
            Given $N = 4$ observed output data points, 
            \begin{equation}
                \vc{y} = \begin{bmatrix} 1 & 4 & 2 & 3\end{bmatrix}^\tr,
            \end{equation}  
            we consider the weighted minimization problem in~\eqref{eq:numericalsystemoptimization} for a first-order model, i.e., $n_a = 1$. 
            We eliminate all variables in the PF system (which consists of twelve equations) except for $s_1$ and $s_2$, leading to one bivariate polynomial equation $t_1(s_1, s_2) = 0$ that contains the Pareto front.
            Its zero set is visualized in \cref{fig:misfitlatency} as a solid blue line (\ref{plot:misfitlatency:pareto}).
            We can verify that the obtained zero set contains the Pareto front by using a sampling approach. 
            Some points of the sampled Pareto front are shown as black dots (~\ref{plot:misfitlatency:sampled}~). 
            Note that, with the approach in \cref{sec:interpretation}, it is now possible to easily retrieve the different weights and model parameters along the Pareto front, providing insights into the models that describe the measured output data.
            
            The two intersection points with the axes correspond to two specific model classes. 
            On the one hand, the intersection with the $s_2$-axis corresponds to a pure autoregressive model without misfit. 
            The value of $s_2$ at this intersection is equal to the norm of the latent input $\norm{\vc{e}}_2^2$. 
            When identifying a pure autoregressive model for the given measured output data (with the globally optimal techniques described in~\cite{vermeersch2019globally}), the same model parameters can be obtained. 
            The corresponding objective values are visualized in \cref{fig:misfitlatency} with the green dot (~\ref{plot:misfitlatency:green}~).
            The intersection with the $s_1$-axis, on the other hand, resembles an autonomous model with output misfit. 
            Computing the misfit in the case that the latent input is zero yields the same value as when using the globally optimal methods from~\cite{lagauw2026exact}, as is shown with the red dot (~\ref{plot:misfitlatency:red}~).
            It is clear that along the Pareto front, the mismatch between the measured output data and model is distributed among these two modifications.
        \end{example}

\section{Conclusions and Future Work}
    \label{sec:conclusion}

    We presented a novel approach to compute the Pareto front in multi-objective minimization for multivariate polynomial optimization problems. 
    Through systematic elimination of all variables except the objective values, the method produces an exact algebraic description that contains the entire Pareto front.
    In contrast to sampling-based methods, our elimination-based approach does not depend on a sampling distance or \textit{a priori} weight assumptions.
    When elimination is performed numerically, it handles floating-point coefficients naturally.
    Several numerical examples confirmed the correctness of the proposed elimination-based approach and illustrate how elimination yields algebraic relations that describe the Pareto front.
    
    We applied elimination also in the context of system identification: computing the Pareto-efficient trade-off between misfit and latency when identifying the model from given output data. 
    This example highlighted an interesting potential application of elimination-based techniques for computing Pareto fronts. 
    The resulting Pareto front confirms the theoretical idea in~\cite{lemmerling2001misfit} that a continuum of equivalent models exists, with the pure latency and pure misfit models as limiting cases. 
    Crucially, the Pareto front gives engineers a principled, complete picture of the trade-off landscape, enabling an informed model selection rather than an arbitrary, \textit{a priori} choice of hyperparameters.
    
    Although our elimination-based approach looks very promising, two challenges that popped up during the system identification problem merit attention in future work. 
    Firstly, the computational complexity remains a bottleneck for problems with many variables or high polynomial degrees; a hierarchical elimination strategy could enable scaling to larger problems. 
    Secondly, extending the elimination idea beyond polynomial objective functions is an important open direction that we intend to pursue in the near future.

    \bibliographystyle{plain}
    \bibliography{bibliography}
\end{document}